\documentstyle{article}
\begin{document}
\baselineskip+8pt
\small \begin{center}{\textit{ In the name of
 Allah, the Beneficent, the Merciful}}\end{center}

\begin{center}{\bf  AN INVERSION FORMULA FOR MULTIVARIATE POWER SERIES}
\end{center}
\begin{center}{\bf  Ural Bekbaev  \\Turin Polytechnic University in Tashkent,\\
INSPEM, Universiti Putra Malaysia.\\
e-mail: bekbaev2011@gmail.com, u.bekbaev@polito.uz}
 \end{center}

\begin{abstract}{For formal multivariate power series $\varphi(x)$ an inversion formula of the form $$ \varphi^{-1}(x)=x +\sum_{m=1}^{\infty}\sum_{k=0}^m (-1)^k\left(\begin{array}{c}
  m \\
  k \\
\end{array}\right)\varphi^{\circ k}(x)\ \ \ \mbox{is offered}$$ }\end{abstract}



\begin{center}{\bf 1. Introduction} \end{center}

In this paper we are going to show that if $\varphi(x): F^n \rightarrow F^n$ is a formal power series of the form $\varphi(x)=x+\mbox{"higher order terms"}$ then
$$ \varphi^{-1}(x)=x +\sum_{m=1}^{\infty}\sum_{k=0}^m (-1)^k\left(\begin{array}{c}
  m \\
  k \\
\end{array}\right)\varphi^{\circ k}(x)$$, where $n$ is any fixed positive integer, $F$ is any field of characteristic zero, $ \varphi^{-1}(x)$ stands for formal inversion of $\varphi$, $\circ$ stands for composition (superposition) operation, $\varphi^{\circ k}(x)$ stands for $k$ times composition of $\varphi$ with itself and $\varphi^{\circ 0}=id$- the identity map $id(x)=x$.

Here we are not going to look for the most general case for which this formula is valid and therefore in future one can assume that $F$ is the field of real or complex numbers. To prove this result we need so called "the symmetric product of matrices" which was introduced before, see for example [1], in a little different form.

Here are the needed results, proofs of which (in polynomial case) can be found in [1].

 For a positive integer $n$ let $I_n$ stand for all row $n$-tuples with nonnegative integer entries with the
 following linear order: $\beta=
 (\beta_1,\beta_2,...,\beta_n)<\alpha=(\alpha_1,\alpha_2,...,\alpha_n)$ if and only if $\vert \beta\vert < \vert \alpha\vert$ or
$\vert \beta\vert = \vert \alpha\vert$ and $\beta_1>\alpha_1$ or $\vert \beta\vert = \vert \alpha\vert$, $\beta_1=\alpha_1$ and $\beta_2>\alpha_2$, et cetera
, where $\vert \alpha\vert$ stands for
$\alpha_1+\alpha_2+...+\alpha_n $. We consider in $I_n$ component-wise addition and subtraction (when the result is in $I_n$), for example, $\alpha +\beta=(\alpha_1+\beta_1,...,\alpha_n+\beta_n)$.  We write $\beta \ll \alpha$ if $\beta_i \leq \alpha_i$
for all $i=1,2,...,n$, $\left(\begin{array}{c}
  \alpha \\
  \beta \\
\end{array}\right)$ stands for $\frac{\alpha!}{\beta!
(\alpha -\beta)!}$, $\alpha !=\alpha_1!\alpha_2!...\alpha_n!$

  For any  nonnegative
integer numbers $p,p'$ let  $M_{n,n}(p',p;F)=M(p',p;F)$ stand for
all $"'p\times p"$ size matrices $A=(A^{\alpha'}_{\alpha})_{\vert
\alpha \vert=p, \vert \alpha'\vert=p'}$ ($\alpha'$ presents
row, $\alpha$ presents column and $\alpha\in I_n,\alpha'\in
I_{n}$). The ordinary size of a such matrix
is $\left(\begin{array}{c}
  p'+n-1 \\
  n-1 \\
\end{array}\right)\times\left(\begin{array}{c}
  p+n-1 \\
  n-1 \\
\end{array}\right)$. Over such kind matrices in addition to the ordinary sum and
product of matrices we consider the following "symmetric product" as well:

{\bf Definition 1.} If $A\in M(p',p;F)$ and $B\in M(q',q;F)$ then
$A\odot B=C\in M(p'+q', p+q;F)$  such that for any
$\vert\alpha\vert=p+q$, $\vert\alpha'\vert=p'+q'$, where
$\alpha\in I_n,\alpha'\in I_{n}$,
$$C^{\alpha'}_{\alpha}=\sum_{\beta,\beta'}\left(\begin{array}{c}
  \alpha \\
  \beta \\
\end{array}\right)
    A^{\beta'}_{\beta}B^{\alpha'-\beta'}_{\alpha-\beta}
$$, where the sum is taken over all $\beta\in I_n,\beta'\in I_{n}$, for which $\vert
\beta\vert=p$, $\vert \beta'\vert=p'$, $\beta\ll \alpha$ and
$\beta'\ll \alpha'$.

%Let us agree that $h$ ($H$, $v$, $V$) stands for any element of
%$M(0,1;F)$ (respect. $M(0,p;F)$, $M(1,0;F)$, $M(p,0;F)$ , where
%$p$ may be any nonnegative integer). The following three propositions
%are about some main properties of the symmetric product, for the proofs one can see [1].

{\bf Proposition 1.} For the above defined product the following
are true.

1. $A\odot B=B\odot A$.

2. $(A+B)\odot C=A\odot C+ B\odot C$.

3. $(A\odot B)\odot C=A\odot (B\odot C)$

4. $ (\lambda A)\odot B=\lambda (A\odot B)$ for any $\lambda\in F$

5.  $A\odot B=0$ if and only if $A=0$ or $B=0$.

In future $A^{\odot m}$ means the $m$-th power of matrix $A$ with
respect to the new product.

{\bf Proposition 2.} If $h=(h_1,h_2,...,h_{n})\in M(0,1;R)$,
 then for any $|\alpha|=m$
$$(h^{\odot m})^0_{\alpha}=m!h^{\alpha}$$, where $h^{\alpha}$ stands for
$h_1^{\alpha_1}h_2^{\alpha_2}...h_n^{\alpha_n}$

{\bf Proposition 3.} For any nonnegative integers $p$, $q$, $p'$,
$q'$ and matrices $A\in M_{n,n}(p',p;F)$, $B\in
M_{n,n}(q',q;F)$, $h=(h_1,h_2,...,h_n)\in M_{n,n}(0,1;F)$,
the following
equality
$$(\frac{h^{\odot p}}{p!}A)\odot
(\frac{h^{\odot q}}{q!}B)=\frac{h^{\odot (p+q)}}{(p+q)!}(A\odot B)
$$ is true.

 In future $Mat_{n,n}(F)=Mat(F)$ stands for the set of all
block matrices $A=(A(p',p))_{p',p}$ with blocks $A(p',p)\in
M_{n,n}(p',p;F)$ for all nonnegative integers $p$, $p'$. In
future it is assumed that $M(p',p;F)$ is a subset of $Mat(F)$ by
identifying each $A(p',p)\in M(p',p;F)$ as the element of $Mat(F)$
which's all blocks are zero, may be, except for $(p',p)$ block
which is $A(p',p)$.

For any $A, B\in Mat(F)$ we  define $A\odot B=C\in Mat(F)$, where
for all nonnegative integers $p$, $p'$
$$C(p',p)=\sum_{q',q}A(q',q)\odot B(p'-q',p-q)$$

The above Propositions show that $(Mat(F); +, \odot)$ is an
integral domain. Its identity element will be $1\in Mat(F)$ whose
all blocks are zero except for $(0,0)$ block which is $1$ -the
identity element of $F$, $F$ is a subring of $Mat(F)$.

 In future the expression
$e^{\odot A}$, whenever it has meaning, stands for
$$1+\frac{1}{1!}A+\frac{1}{2!}A^{\odot 2}+\frac{1}{3!}A^{\odot 3}+ . .
.=\sum_{i=0}^{\infty}\frac{1}{i!}A^{\odot i}$$, $F[[x]]$ is the
ring of formal power series in variables $x_1,x_2,...,x_n$  over
$F$, $x=(x_1,x_2,...,x_n)\in M_{n,n}(0,1;F[[x]])$.

If $\varphi(x)=(\varphi_1(x),\varphi_2(x),...,\varphi_{n}(x))\in
F[[x]]$
 one can screen it in the form
$$\varphi(x)=x^{\odot 0}M^0_1+\frac{x^{\odot 1}}{1!}M_1^1+
\frac{x^{\odot 2}}{2!}M_1^2+...=\sum_{i=0}^{\infty}\frac{x^{\odot i}}{i!}M^i_1=
e^{\odot x}M_{\varphi} $$, where $x^{\odot 0}=1$, $M^{p'}_1\in Mat(p',1;F)$,
$M_{\varphi}\in Mat(F)$ with blocks $M_{\varphi}(p,p')$ such that
$M_{\varphi}(p',p)=0$ whenever $p\neq 1$ and
$M_{\varphi}(p',1)=M^{p'}_1$ for all nonnegative integers. We call
$M_{\varphi}$ the matrix of $\varphi(x)$.

Let us consider only power series with zero constant terms.

The following theorem deals with the matrix of composition of
power series.

{\bf Theorem 1.} If
$\psi(x)=(\psi_1(x),\psi_2(x),...,\psi_n(x))=e^{\odot
x}M_{\psi}\in F[[x]]$ then $$M_{\psi \circ \varphi}=e^{\odot
M_{\varphi}} M_{\psi}$$

The associative property of composition yields in the following
result.

 {\bf Theorem 2.} If $\xi(x)=(\xi_1(x),\xi_2(x),...,\xi_{n}(x))=e^{\odot x}M_{\xi} \in
F[[x]]$ then
$$e^{\odot M_{\varphi}}(e^{\odot M_{\psi}}M_{\xi})=e^{\odot (e^{\odot M_{\varphi}}M_{\psi})}M_{\xi}$$

{\bf Corollary 1.} For any natural $m$ one has
$$M_{\varphi^{\circ m}}=(e^{\odot M_{\varphi}})^{m-1}M_{\varphi}=(e^{\odot M_{\varphi}})^{m}E_1$$,  where $E_1$ is $"1\times 1"$ size identity
matrix.

\begin{center}{\bf 2. The main result} \end{center}

In future we consider any fixed $\varphi(x)=e^{\odot
x}M_{\varphi}$ for which $M_{\varphi}(0,1)=0$,
$M_{\varphi}(1,1)=E_1$, where $M_{\varphi}(p',1)=M^{p'}_1\in Mat(p',1;F)$ are arbitrary for $p'\geq
2$. Let $\varphi^{-1}(x)=e^{\odot x}M_{\varphi^{-1}}$ stand for
inverse formal power series to $\varphi$. The block components of
$M_{\varphi^{-1}}$ we denote by $N^{p'}_1$

{\bf Theorem 3.} For the above mentioned power series $\varphi(x)$
the following equality
$$ \varphi^{-1}(x)=\sum_{m=0}^{\infty}\sum_{k=0}^m (-1)^k\left(\begin{array}{c}
  m \\
  k \\
\end{array}\right)\varphi^{\circ k}(x)=x +\sum_{m=1}^{\infty}\sum_{k=0}^m (-1)^k\left(\begin{array}{c}
  m \\
  k \\
\end{array}\right)\varphi^{\circ k}(x)$$ is true.

{\bf Proof.} Due to Theorem 1 equality $\varphi^{-1}(\varphi(x))=x
$  is nothing than
 $e^{\odot M_{\varphi}}M_{\varphi^{-1}}=E_1$ and therefore $ M_{\varphi^{-1}}=(e^{\odot M_{\varphi}})^{-1}E_1$ provided that $e^{\odot M_{\varphi}}$
 is invertible.

 But  in our case $$(e^{\odot M_{\varphi}})^{-1}= (E_{\infty}- (E_{\infty}-e^{\odot M_{\varphi}}))^{-1}
 =E_{\infty}+ \sum_{m=1}^{\infty}(E_{\infty}-e^{\odot M_{\varphi}})^m$$ is well defined, where $E_{\infty}$ stands for infinite size identity matrix.

%Note that $(E_{\infty}-e^{\odot M_{\varphi}})E_1=E_1-M_{\varphi}=
%M_{id-\varphi }$

Due to Corollary 1 one can see that $$
(E_{\infty}-e^{\odot M_{\varphi}})^mE_1=\sum_{k=0}^{m}(-1)^{k}\left(\begin{array}{c}
  m \\
  k\\
\end{array}\right)M_{\varphi^{\circ k}}$$, where $\varphi^{\circ 0}=id$ -the
identity map $id(x)=x$.

It implies that
$$ M_{\varphi^{-1}}=E_1+ \sum_{m=1}^{\infty}\sum_{k=0}^{m}(-1)^{k}\left(\begin{array}{c}
  m \\
  k\\
\end{array}\right)M_{\varphi^{\circ k}}$$

Now one can write this result in terms of $\varphi$
$$\varphi^{-1}(x)=x+ \sum_{m=1}^{\infty}\sum_{k=0}^{m}(-1)^{k}\left(\begin{array}{c}
  m \\
  k\\
\end{array}\right)\varphi^{\circ k}(x)$$ or in a symbolic form $$\varphi^{-1}(x)=x+ \sum_{m=1}^{\infty}(id-\varphi)^{[\circ]m}(x) $$,
where  $ (id-\varphi)^{[\circ]m}(x)$ stands for
$\sum_{k=0}^{m}(-1)^{k}\left(\begin{array}{c}
  m \\
  k\\
\end{array}\right)\varphi^{\circ k}(x)$, that is one can remove parentheses in $ (id-\varphi)^{[\circ]m}(x)$ as if $\varphi$ were a linear operator.
This is the proof of Theorem 3.

One can ask the following question.

{\bf Question.} If both of $\varphi(x)$, $\varphi^{-1}(x)$ are
polynomial maps does it imply that for some $m_0> 0$
$$\sum_{m=m_0}^{\infty}\sum_{k=0}^{m}(-1)^{k}\left(\begin{array}{c}
  m \\
  k\\
\end{array}\right)\varphi^{\circ k}(x)=0 \ ?$$

In common case I am not sure that the answer to this question is positive.

Due to $$\sum_{k=0}^{m+1}(-1)^{k}\left(\begin{array}{c}
  m+1 \\
 k\\
\end{array}\right)(e^{\odot M_{\varphi}})^{k}E_1=(E_{\infty}-e^{\odot M_{\varphi}})^{m+1}E_1=(E_{\infty}-e^{\odot M_{\varphi}})^{m}M_{id-\varphi}=$$
$$\sum_{k=0}^{m}(-1)^{k}\left(\begin{array}{c}
  m \\
 k\\
\end{array}\right)(e^{\odot M_{\varphi}})^{k}M_{id-\varphi}=
\sum_{k=0}^{m}(-1)^{k}\left(\begin{array}{c}
  m \\
 k\\
\end{array}\right)M_{(id-\varphi)\circ \varphi^{\circ k}}=\sum_{k=0}^{m}(-1)^{k}\left(\begin{array}{c}
  m \\
  k\\
\end{array}\right)(M_{\varphi^{\circ k}}-M_{\varphi^{\circ (k+1)}})$$  for
$\Phi_m(x)=\sum_{k=0}^{m}(-1)^{k}\left(\begin{array}{c}
  m \\
  k\\
\end{array}\right)\varphi^{\circ k}(x)$ one has
$\Phi_{m+1}(x)=\Phi_m(x)-\Phi_m(\varphi(x))$, where $\Phi_0 =id$.

{\bf Proposition 4.} For any $m_0\geq 0$ the following equality  $$\Phi_{m_0}(x)=\sum_{m=m_0}^{\infty}(\Phi_{m}\circ \varphi)(x)= \sum_{m=m_0}^{\infty}\Phi_{m}(\varphi(x))$$ is true.

{\bf Proof.} At $m_0 > 0$ due to the equality
$\Phi_{m+1}(x)=\Phi_m(x)-\Phi_m(\varphi(x))$ one has $$ \sum_{m=m_0}^{\infty}\Phi_{m}(x)=\sum_{m=m_0}^{\infty}(\Phi_{m-1}(x)-\Phi_{m-1}(\varphi(x)))=$$
$$\Phi_{m_0-1}(x)-\Phi_{m_0-1}(\varphi(x))+\sum_{m=m_0+1}^{\infty}(\Phi_{m-1}(x)-\Phi_{m-1}(\varphi(x)))=\Phi_{m_0}(x)+ \sum_{m=m_0}^{\infty}\Phi_{m}(x)-\sum_{m=m_0}^{\infty}\Phi_{m}(\varphi(x))$$
, which implies that $\Phi_{m_0}(x)=\sum_{m=m_0}^{\infty}\Phi_{m}(\varphi(x))$. The case $m_0=0$ is a consequence of the case $m_0=1$.

 {\bf Corollary 2.} If $\sum_{m=m_0}^{\infty}\Phi_{m}(x)=0$ for some $m_0\geq 1$ then $\Phi_{m}(x)=0$
 for any $m\geq m_0$.

 {\bf Corollary 3.} If $\varphi(x)$ is a polynomial map
and for some $m\geq 1$ the equality
$\sum_{k=0}^{m}(-1)^{k}\left(\begin{array}{c}
  m \\
  k\\
\end{array}\right)\varphi^{\circ k}(x)=0$ is true then  $\varphi^{-1}(x)$ is also a polynomial map.

In an equivalent form the condition of Corollary 3 can be given
in the form: For some $m\geq 1$ the equality $$(\partial
\varphi)(x)(mE_1+\sum_{k=2}^{m}(-1)^{k-1}\left(\begin{array}{c}
  m \\
  k\\
\end{array}\right) (\partial
\varphi)|_{\varphi(x)}(\partial \varphi)|_{\varphi^{\odot
2}(x)}...(\partial \varphi)|_{\varphi^{\circ (k-1)}(x)})=E_1$$
is true, where $\partial_k=\frac{\partial}{\partial x_k}$,
$(\partial\varphi(x))^i_j=\partial_i\varphi_j(x)$, $(\partial \varphi)|_{\varphi^{\circ (k-1)}(x)}=(\partial \varphi)(\varphi^{\circ (k-1)}(x))$.

 For block components of $M_{\varphi^{-1}}$ the following recurrent formula is true.

{\bf Proposition 5.} For any $m>1$ one has
$$N^m_1=-\sum_{k=1}^{m-1}(\sum_{|\alpha|=k,
\|\alpha\|=m}\frac{(M_1^1)^{\odot \alpha_1}\odot (M^2_1)^{\odot
\alpha_2}\odot (M^3_1)^{\odot \alpha_3}\odot ...}{\alpha!}
)N^k_1$$ and $N^1_1=M^1_1$, where
$\|\alpha\|=1\alpha_1+2\alpha_2+3\alpha_3+...$.

{\bf Proof.} The above equalities are nothing than the equality $e^{\odot M_{\varphi}}M_{\varphi^{-1}}=E_1$ in $(m,1)$ blocks for $m>1$.

{\bf Remark.} The set of polynomial maps $\varphi(x)$ for which
$\varphi^{-1}(x)$ is also polynomial map is a group with respect
to the composition operation. Finding any system of generators of it may be useful.

\begin{center}{References}\end{center}

[1] Ural Bekbaev. \emph{A matrix representation of composition of
polynomial maps}. \\ arXiv0901.3179v3 [math. AC] 22 Sep.2009

\end{document}